\documentclass[a4paper]{iopconfs}
\usepackage{graphicx}
\usepackage{amsmath}
\usepackage{epsfig}
\usepackage{amssymb}

\def\varphi0{{\cal Q}}

\def\ox{\bar{x}}
\def\oX{\bar{X}}
\def\ot{\bar{t}}
\def\ol{\bar{\lambda}}

\renewcommand{\thesubsection}{\arabic{section}.\arabic{subsection}}
\renewcommand{\thesection}{\arabic{section}}

\newcommand{\SE}{\setcounter{equation}{0} \section}
\newcommand{\be}{\begin{equation}}
\newcommand{\ee}{\end{equation}}
\newcommand{\baa}{\begin{array}}
\newcommand{\eaa}{\end{array}}
\newcommand{\ba}{\begin{eqnarray}}
\newcommand{\ea}{\end{eqnarray}}

\newtheorem{theo}{\textsc{Theorem}}[section]

\begin{document}
\title{Homogenization of dislocation dynamics}

\author{Ahmad El Hajj, Hassan Ibrahim and R\'egis Monneau}

\address{CERMICS, ENPC, 6 \& 8 avenue Blaise
    Pascal, Cit\'e Descartes, Champs sur Marne, 77455 Marne-la-Vall\'ee
    Cedex 2, France}

\ead{elhajj@cermics.enpc.fr, ibrahim@cermics.enpc.fr, monneau@cermics.enpc.fr}

\begin{abstract}
In this paper we consider the dynamics of dislocations with the same
Burgers vector, contained in the same glide plane, and moving in a material
with periodic obstacles. We study two cases: i) the
particular case of parallel straight dislocations and ii) the general case of curved
dislocations. In each case, we perform rigorously the homogenization of the
dynamics and predict the corresponding effective macroscopic elasto-visco-plastic flow rule.
\end{abstract}



\SE{Introduction}

In the recent years, an important effort has been done, both to improve the
methods to compute discrete dislocation dynamics (see for instance the book
of Bulatov and Cai
\cite{BC} and the references therein) and also to connect them
to continuum models of plasticity in crystalline solids (see for instance
Fivel et al. \cite{FTRC} and more recently Hoc et
al. \cite{HDK}). Although continuum models of 
dislocations are known since the 50's (see Kr\"{o}ner \cite{K,K2}), the
dynamics has been taken into account  only recently : see Groma et
al. \cite{GB,GCZ} in 2D (and their mathematical studies in
\cite{FE,IJM}), Hochrainer et al. \cite{HZG}, and Monneau  \cite{M} in 3D. The goal
of our work is to present, on a particular 
example,  a rigorous justification of a continuum model with densities of
dislocations  bridging the
gap with dislocation dynamics at the microscale. Indeed for a
very special geometry, we are able to deduce by homogenization, the
macroscopic elasto-visco-plastic flow rule relating the plastic strain
velocity to the shear stress. The full technical details are presented
in \cite{FIM}.

\SE{Homogenization of straight dislocations}

In this section, we consider the case of parallel straight edge dislocations
with the same Burgers vector $\mbox{\bf b}=be_x$ with $b>0$, where
$(e_x,e_y,e_z)$ is an orthonormal basis with corresponding coordinates
$(x,y,z)$. All these dislocation lines are 
assumed to be contained in the same glide plane $(x,y)$ and to move in
this plane.

\subsection{The microscopic model for straight dislocations}

Because of our assumptions, for every integer $i\in \mathbb{Z}$, we can
simply describe the position of the $i$-th dislocation
by its real abscissa that we call $x_i(t)$ where $t$ is the time. We want
to take into account the interactions of each dislocation with other
defects in the crystal, that constitute obstacles
to their motion. Those obstacles can be for instance other pinned
dislocations or precipitates. In order to simplify the analysis, we will
assume that these obstacles are periodically distributed, of spatial period
$\lambda$. In our model, those obstacles will be simply modeled by a smooth
periodic potential $V^{per}$ satisfying
$V^{per}(x+\lambda)=V^{per}(x)$. Then the energy of the system is the sum
of two contributions: the interactions of each dislocations with the
periodic potential and the sum of the two-body interactions between
dislocations associated to a pair potential $V$.
The energy of a set of
dislocations is then
given by
$$E=\sum_{i} V^{per}(x_i) + \sum_{i<j} V(x_i-x_j)\quad  
\mbox{with}\quad V(x)=-\bar{\mu} b\ln |x| \quad \mbox{and}\quad
\bar{\mu}=\frac{\mu }{2\pi (1-\nu)}$$ 
where the constants $\mu$ and $\nu$ are respectively the shear modulus
and the Poisson ratio. Remark that 
the force $-V'(x)$ is then the usual Peach-Koehler force created at the
point $x$ by an edge dislocation positioned at the origin. 

We then consider the fully overdamped dynamics, where the velocity is
proportional to the force, i.e.
\begin{equation}\label{eq::1}
B\frac{dx_i}{dt}=-\nabla_{x_i} E + \tau_{ext}   
\end{equation}
where $B$ is the viscous drag coefficient and the force is 
on the right hand side. The first contribution to the force is a term  deriving from the
energy and $\tau_{ext}$ is a real exterior applied shear stress, that can be
seen as a driving force of the system. 
A natural question is then: what is the macroscopic behavior of this system ?

In order to answer this question (which is done in Theorem \ref{th::1}),
we have to introduce the plastic strain. To each dislocation is associated a three-dimensional
displacement in the crystal, whose plastic strain is localized in the glide
plane $z=0$ and is equal to $\gamma \delta_0(z)$ where $\delta_0$ is the
Dirac mass. For instance, for a dislocation $x_i$, the intensity
$\gamma$ (that we continue to call plastic strain) is equal to $-b
H(x-x_i)$ where the Heaviside 
function $H(x)$ is equal to $1$ for positive $x$ and zero otherwise. Here 
the sign defining the plastic strain is such that the quantity
$\gamma$ increases when $x_i$ increases. Then the total plastic strain can
be written as
$$\gamma(x,t)=-b \sum_{i} H(x-x_i(t)).$$

\subsection{The normalization procedure}

We are now interested in the behavior of the system at a macroscopic scale
$\Lambda$ such that $\Lambda >> \lambda =\ol b$ where $\ol >1$ is a fixed
ratio. Then we introduce several dimensionless quantities. We call
$\bar{x}$ and $\ot$ the normalized
spatial and time coordinates at the macroscopic level, and introduce a
parameter $\varepsilon$ and the associated
normalized macroscopic plastic strain
${\gamma}^\varepsilon$ such that
\begin{equation}\label{eq::0}
\bar{x} = \frac{x}{\Lambda}, \quad \ot =
\frac{\bar{\mu}}{B} \frac{t}{\Lambda},\quad \varepsilon = \frac{b}{\Lambda}  \quad \mbox{and}\quad 
\displaystyle{{\gamma}^\varepsilon(\ox,\ot)=
  \frac{1}{\Lambda}\gamma(x,t)} \quad \mbox{with}\quad {\gamma}^\varepsilon(\ox,0)=\varepsilon\left[\frac{1}{\varepsilon}\gamma_0(\ox)\right]   
\end{equation}
where $\left[\cdot\right]$ is the floor function, $\gamma_0$ is a given function  and $B\Lambda/\bar{\mu}$ is a
typical macroscopic time deduced from equation (\ref{eq::1}). 
Remark that $\varepsilon$ can be very small in our application (for instance $\varepsilon \simeq 10^{-6}$ if
$b\simeq 10^{-9}m$ and $\Lambda\simeq 10^{-3}m$).

We expect that the macroscopic behavior of the model is well described by
{\it the limit macroscopic plastic strain $\gamma^0(\ox,\ot)$ of $\gamma^\varepsilon(\ox,\ot)$ as
  $\varepsilon$ goes to zero}.

\subsection{Heuristics for the macroscopic stress field}

  In this subsection, we want to give heuristic expressions of the normalized
  dislocation density and the macroscopic stress field, in terms of the limit
  macroscopic plastic strain.

 We remark that the gradient of the map $x\mapsto -\gamma^\varepsilon
  (x/\Lambda,\ot) /\varepsilon$
 is a sum of Dirac masses, and then the number of dislocations in a large segment
 of length $\Delta x$ is formally given by $-\int_{0}^{\Delta x} 
 \frac{1}{\varepsilon \Lambda } \frac{\partial \gamma^\varepsilon}{\partial \ox}(x/\Lambda,\ot)\
 dx$. This shows at least formally that the dislocation density can be
 estimated as $\rho(x,t)=-\frac{1}{\varepsilon\Lambda} \frac{\partial \gamma^0}{\partial
   \ox}(\ox,\ot)$.
 Then the total stress on the right hand side of (\ref{eq::1}) can be formally
 described at the macroscopic scale by
  \begin{equation}\label{eq::2}
 \displaystyle{\tau = \tau_{ext} + \tau_{sc} \quad \mbox{with}\quad 
  \tau_{sc}(\ox,\ot)= -\bar{\mu} \int_{-\infty}^{+\infty}
   \frac{d\ox'}{\ox-\ox'} \frac{\partial \gamma^0}{\partial \ox}(\ox',\ot)}
  \end{equation}
 where we take the principal value in the integral defining the
 self-consistent field $\tau_{sc}$. This expression can be deduced from the
 equation $\tau_{sc}(\ox,\ot)= -(V'\star_x \rho)(x,t)$, where $\star_x$ denotes the
 convolution with respect to the variable $x$.
Remark also that the expression (\ref{eq::2}) of
 $\tau_{sc}$ 
is known to be the resolved shear stress created by the normalized
 dislocation density
 \begin{equation}\label{eq::3}
 \rho^0=- \frac{\partial \gamma^0}{\partial \ox}   
 \end{equation}
where for instance $\rho^0=1/\bar{\lambda}$ when there is one dislocation by spatial period $\lambda$.
In particular, we see that $\tau_{sc}$ keeps  the memory of the long range interactions
 between dislocations.

\subsection{The homogenization result}\label{s1.3}


We expect that the effective equation satisfied by the limit
$\gamma^0$ can be written
\begin{equation}\label{eq::5}
\left\{\begin{array}{l}
\displaystyle{\frac{\partial \gamma^0}{\partial \ot} = f(\rho^0, \tau), \quad
\mbox{for all}\quad \ox\in\mathbb{R},\quad \ot\in (0,+\infty)},\\
\\
\gamma^0(\ox,0)=\gamma_0(\ox) \quad \mbox{for all}\quad \ox\in \mathbb{R}
\end{array}\right. 
\end{equation}
where $\rho^0$ is given in (\ref{eq::3}) and $\tau$ in (\ref{eq::2}). 
Then our main result is:
\begin{theo}\label{th::1}{\bf (Homogenization of straight dislocations)}\\
Assume that the initial data $\gamma_0$ is non-decreasing and satisfies
$|\gamma_0|+ |\gamma_0'|+ |\gamma_0''| \le C$ for some constant $C$.
Then for any $C^2$ periodic potential $V^{per}$, there exists
a continuous function $f: \mathbb{R}^2 \to \mathbb{R}$ such that $\tau \mapsto
f(\rho^0,\tau)$ is nondecreasing. And there exists a unique viscosity
solution $\gamma^0$ of the equation (\ref{eq::5}).\\
Moreover, under the assumptions and notation of this
section, there exists a unique solution $\gamma^\varepsilon$ associated to
the dynamics  (\ref{eq::1})  with initial
data given in (\ref{eq::0}), and $\gamma^\varepsilon$ converges to $\gamma^0$ locally
uniformly on $\mathbb{R} \times [0,+\infty )$.
\end{theo}

This result is proven rigorously in \cite{FIM} in the mathematical framework of
viscosity solutions (see for instance Crandall, Ishii, Lions \cite{CIL} for an introduction to
this theory). We explain in the next section how we compute the
function $f$, which keeps the memory of the short range interactions
between the dislocations and the periodic potential $V^{per}$.

\subsection{Computation of $f$ using Orowan's law}\label{ssf}

In this subsection, we briefly explain (without any justifications) how to
compute the function $f$. We refer the reader to \cite{FIM} for the proofs
of those results.\\
\noindent {\it Case A:} $V^{per}\equiv 0$.\\
In this special case, we can show that
\begin{equation}\label{eq::6}
f(\rho^0,\tau)=\rho^0 \bar{v} \quad \mbox{with}\quad \bar{v}=\frac{\tau}{\bar{\mu}}   
\end{equation}
which is nothing else than the normalized Orowan's law  giving, in a
dimensionless form, the
plastic strain velocity as the product of the normalized dislocation
density $\rho^0$  and the normalized mean velocity $\bar{v}$ of the
dislocations.\\

\noindent {\it Case B:} General periodic potential $V^{per}$.\\
In that case, the function $f$ can
be computed using the following two steps.\\
\noindent {\underline{Step 1}.}\\
For $i\in\mathbb{Z}$, we look for solutions to (\ref{eq::1}) of the following
special form
$$x_i(t)=b\cdot h\left(\frac{vt}{b} + \frac{i}{\rho^0}\right), \quad \mbox{with}\quad
h(a+\bar{\lambda})=\bar{\lambda}+h(a) \quad \mbox{for all}\quad a\in\mathbb{R}$$
for some constant $v$ and for a function $h$ which is called a {\it hull
  function}. Both $v$ and $h$ have to be determined.
Because of the convexity of the two-body potential $V$ outside the origin, 
it is possible to show  that the constant $v$ exists and is unique. Moreover
this constant $v$ can be interpreted as the mean velocity of each dislocation.\\
\noindent {\underline{Step 2}.}\\
We simply define $f(\rho^0,\tau_{ext})$ using the normalized Orowan's law as in
(\ref{eq::6}), but with the normalized velocity
$\bar{v}$ replaced by the constant $\bar{v}=\frac{B}{\bar{\mu}}v$.

\subsection{Numerical computation of $f$}

\begin{figure}[!h]
     \begin{minipage}[b]{.46\linewidth}
      \centering\epsfig{figure=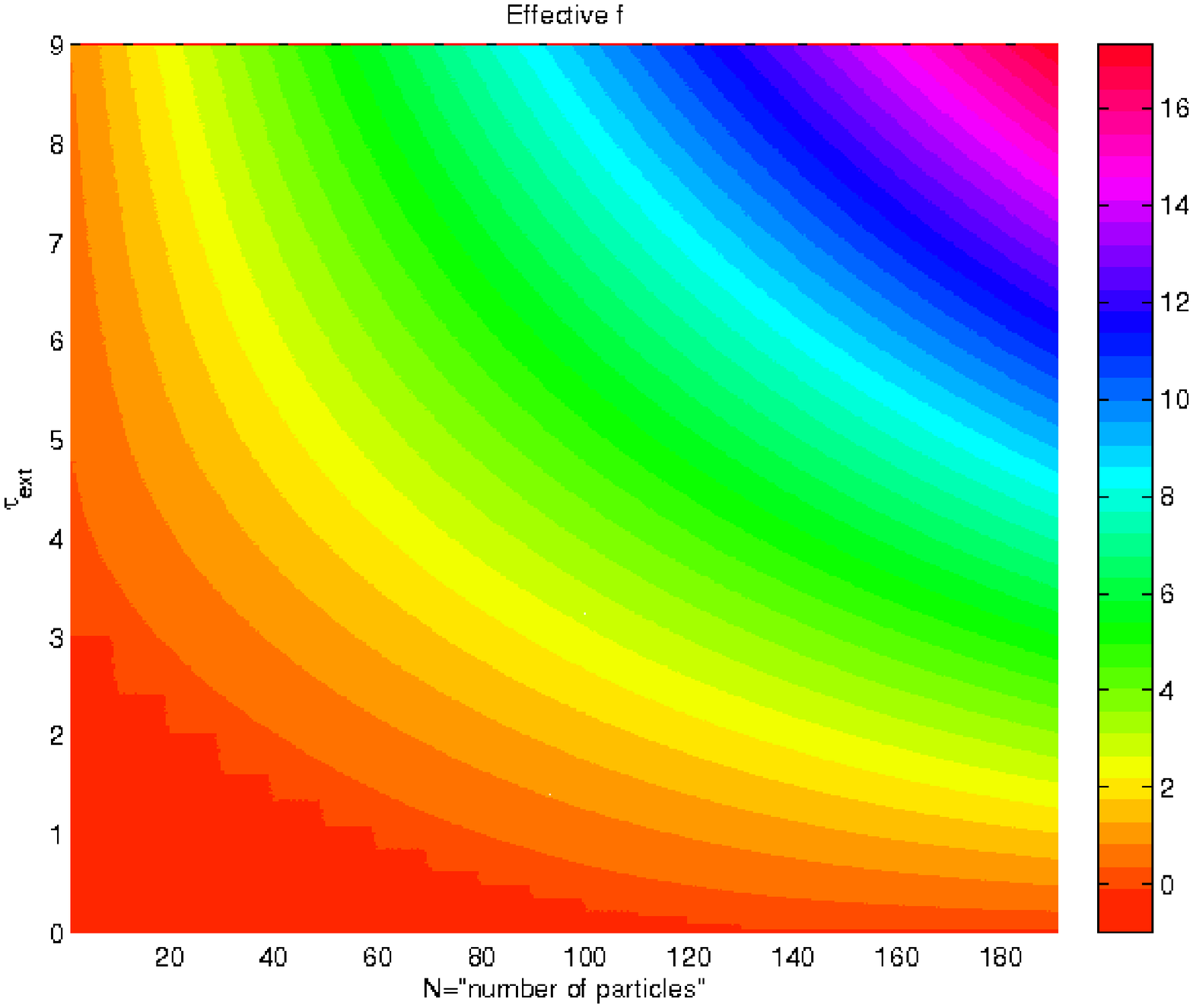,width=\linewidth}
      \caption{Level sets of the effective $f(N/l,\tau_{ext})$ with $N$ on
        abscissas and  $\tau_{ext}$ on ordinates\label{F1}}
     \end{minipage} \hfill
     \begin{minipage}[b]{.46\linewidth}
      \centering\epsfig{figure=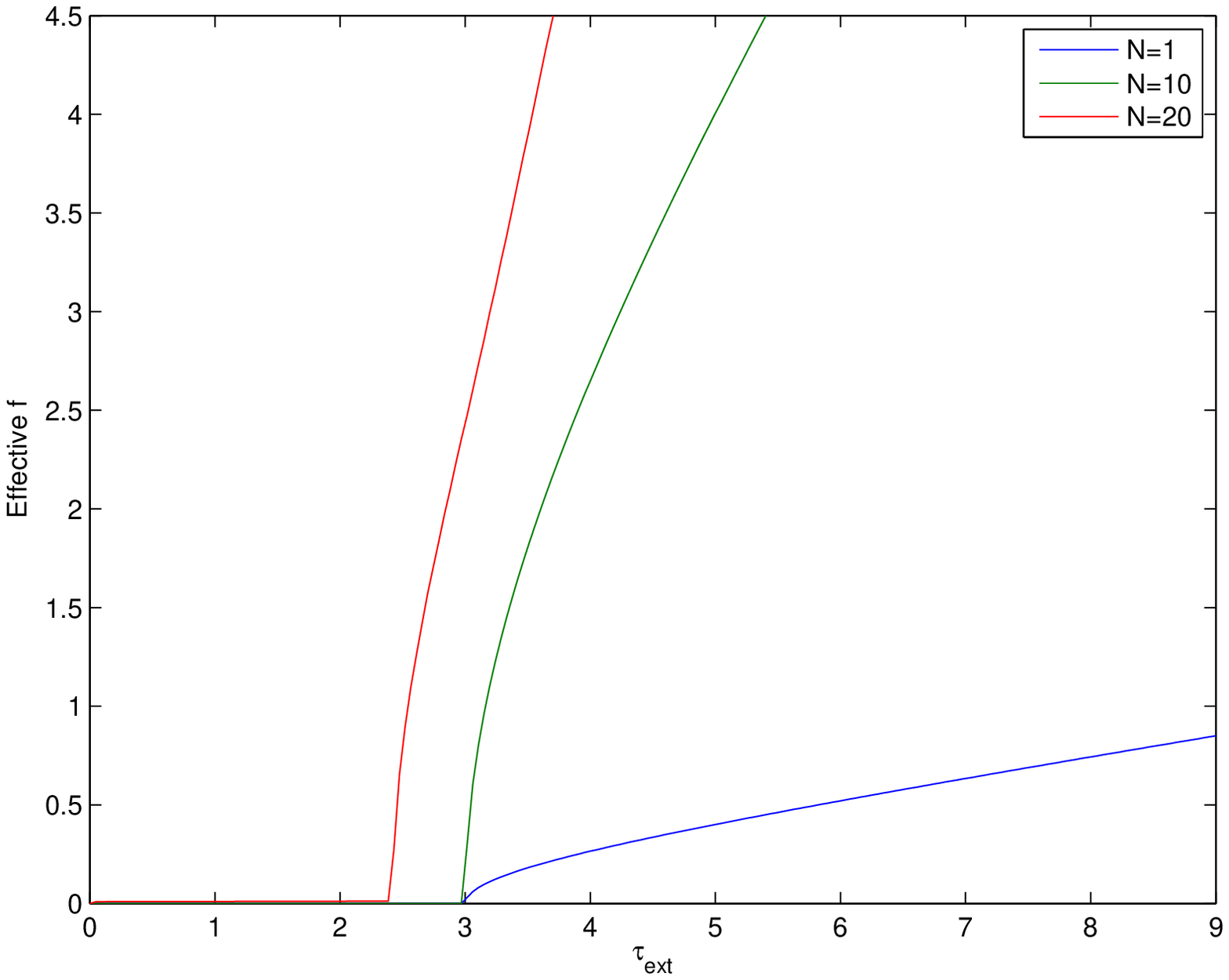 ,width=\linewidth}
      \caption{For  $N=1,10,20$, graph of the map $\tau_{ext}\mapsto f(N/l,\tau_{ext})$\\\label{F2}}
     \end{minipage}
    \end{figure}

We present numerical simulations for the computation of the function
$f$. We work with dimensionless quantities: 
$\lambda=1=\bar{\lambda}=b=B=\bar{\mu}$. 
We put initially $N$ dislocations in
an interval of length $l=10$ which is repeated periodically. Therefore this
interval contains $l$ times the period of
the periodic potential that we choose equal to $V^{per}(x)=\frac{A}{2\pi}\sin(2\pi x)$
with $A=3$. We discretize the ODE system (\ref{eq::1}), using an
explicit Euler scheme with a time step  $\Delta t =0.01$. We compute
numerically the mean velocity $v$ of the dislocations  after a final time
$T=1000$. We then set $f=\rho^0 v$ with $\rho^0=N/l$. We do the computation
with $N=1,...,200$ and $0\le \tau_{ext}\le 9$ with $\Delta \tau_{ext} =
\frac{9}{200}$. Remark that
we can restrict our computation for positive $\tau_{ext}$, because we have
$f(\rho^0,-\tau_{ext})=-f(\rho^0,\tau_{ext})$, from the symmetry of the
potential $V^{per}$ in our problem. The
level sets of the function $f$ are represented on Figure \ref{F1}. In
order to have a better view of the set where $f=0$, this set is
conventionally represented in Figure \ref{F1} with artificial negative
values of $f$. We remark that this figure shows in particular a
collective behavior of the dislocations: higher is the density of
dislocations, then easier the dislocations move above the obstacles.

Figure \ref{F2} shows the map $\tau_{ext}\mapsto f(\rho^0,\tau_{ext})$ for
$\rho^0=N/l$ with $N=1,10,20$. We see in particular that for $\tau_{ext}$ under a threshold
(that depends on the dislocation density $\rho^0$) the function $f$ vanishes.

 \SE{Homogenization of curved dislocations}

 In this section, we very briefly generalize the previous analysis to the
 case of  curved dislocations all contained in the same plane $(x,y)$ with the
 same  Burgers vector $\mbox{\bf b}=be_x$ with $b>0$.

 \subsection{The microscopic model for curved dislocations}

 For $i\in\mathbb{Z}$, the motion of the $i$-th dislocation curve $\Gamma_i(t)$ at the point
 $X\in\mathbb{R}^2$ is given by its
 normal velocity ${\mathcal V}$ defined by
 \begin{equation}\label{eq::1bis}
 B\cdot{\mathcal V}(X,t)= \tau^{per}(X) + \sum_{j} F_{j}(X,t)   
 \end{equation}
 where $F_{j}(X,t)$ is the resolved Peach-Koehler force created by the
 dislocation $\Gamma_j(t)$ at the point $X$. Here $\tau^{per}$ is a smooth
 periodic function satisfying
 $\tau^{per}(X +\lambda k)=\tau^{per}(X)$ for all $k\in\mathbb{Z}^2$, which
 represents the periodic obstacles to the motion of the dislocations and can also
 include the exterior applied stress. To give the expression of this
 force, it is convenient 
 to introduce a continuous function $\tilde{\gamma}(X,t)$
 such that each dislocation curve $\Gamma_j(t)$ can be seen as the level set
 $\tilde{\gamma}(X,t)=jb$ (when this level set is non-degenerated). Then a
 good approximation is given by
 $$\displaystyle{F_{j}(X,t)=\frac12 \int_{\mathbb{R}^2} dZ\ J(X-Z)\ \mbox{sign}(\tilde{\gamma}(Z,t)-jb)}$$
 where, in the integral, the sign function takes values $-1,0,1$. Here the
 kernel $J$ is smooth and satisfies for a cut-off radius $R=\bar{R}b$ with $\bar{R}>1$ fixed:
 $$J(-X)=J(X)\ge 0,  \quad  \mbox{and}\quad 
 \displaystyle{J(X)=J_{\infty}
 (X):=\frac{1}{|X|^3}g\left(\frac{X}{|X|}\right) \quad \mbox{for}\quad
 |X|>R>0}$$
where for isotropic elasticity with $X=(x,y)$, we have 
$g\left(\frac{X}{|X|}\right)=\frac{\mu b}{4\pi}\left\{\frac{x^2(2\beta -1) +
 y^2(2-\beta)}{x^2+y^2}\right\}$ with $\beta=\frac{1}{1-\nu}$. Remark that
 this formula allows to describe with the same formalism edge, screw and
 mixed dislocations (see for instance \cite{AHLM}).
We also define the plastic strain $\gamma$ as 
 $$\gamma= b \left[
   \frac{\tilde{\gamma}}{b}\right]$$ 
 where we recall that $\left[\cdot\right]$ is the
 floor function. Then we proceed as in the previous section and define
 \begin{equation}\label{eq::0bis}
 \bar{X}=\frac{X}{\Lambda},\quad 
 \bar{t}=\frac{\mu}{B}\frac{t}{\Lambda},\quad
 \varepsilon=\frac{b}{\Lambda},\quad \mbox{and}\quad 
 \gamma^{\varepsilon}(\bar{X},\ot)=\frac{1}{\Lambda}\gamma(X,t),\quad  \mbox{with}\quad
 \gamma^{\varepsilon}(\bar{X},0)=\varepsilon
 \left[\frac{1}{\varepsilon}\gamma_0(\oX)\right].  
 \end{equation}

 \subsection{The homogenization result}

 We expect that the effective equation satisfied by the limit
 $\gamma^0$ of $\gamma^\varepsilon$ can be written
 \begin{equation}\label{eq::5bis}
 \left\{\begin{array}{l}
 \displaystyle{\frac{\partial \gamma^0}{\partial \ot} = f(-\nabla \gamma^0, \tau_{sc}), \quad
 \mbox{for all}\quad \oX\in\mathbb{R}^2,\quad \ot\in (0,+\infty)},\\
 \\
 \gamma^0(\oX,0)=\gamma_0(\oX) \quad \mbox{for all}\quad \oX\in \mathbb{R}^2
 \end{array}\right. 
 \end{equation}
 with
 $$\tau_{sc}(\oX,\ot)=\int_{\mathbb{R}^2} dZ\ J_\infty (\oX-Z) \gamma^0(Z,\ot)$$
 where we take the principal value of the integral. Remark that this
 expression of $\tau_{sc}$ is consistent with the one given in (\ref{eq::2})
 in the special case where $\gamma^{0}(\bar{x}, \bar{y}, \bar{t})$ is
 independent of $\bar{y}$.
 Then we have
 \begin{theo}\label{th::1bis}{\bf (Homogenization of curved dislocations)}\\
 Assume that the initial data satisfies $|\gamma_0|+ |\nabla \gamma_0|+
 |D^2 \gamma_0|\le C$ for some constant $C$.
Then for any $C^2$ periodic  function $\tau^{per}$, there exists
 a continuous function $f: \mathbb{R}^2\times \mathbb{R} \to \mathbb{R}$ such that $\tau \mapsto
 f(\cdot,\tau)$ is nondecreasing. And there exists a unique viscosity
 solution $\gamma^0$ of the equation (\ref{eq::5bis}).\\
 Moreover, under the assumptions and notation of this
 section, there exists a unique solution $\gamma^\varepsilon$ associated to
 the dynamics  (\ref{eq::1bis})  with initial
 data given in (\ref{eq::0bis}), and $\gamma^\varepsilon$ converges to $\gamma^0$ locally
 uniformly on $\mathbb{R}^2 \times [0,+\infty )$.   
 \end{theo}

 \SE{Conclusion}

 The main result of our work is the justification of the
 elasto-visco-plastic flow rule by the homogenization of the
 dynamics of dislocations with the same
 Burgers vector, moving in the same glide plane with periodic obstacles. 
 Even if this geometry is very particular, this is, up
 to our knowledge, the first rigorous result in this direction. We also
 explained how to compute the flow rule, and presented numerical results.
 The proof of the homogenization for straight dislocations uses strongly the local convexity of the
 two-body potential $V$  (which is equivalent to the non-negativity of the
 kernel $J$ in the case of curved dislocations). 

 Remark that for the same dynamics, it is
 possible to find non-convex potentials $V$, for which there is no
 homogenization. For a general geometry, there is in general no hope to find any
 convexity argument to justify homogenization. On the contrary, it seems
 reasonable to think that homogenization could arise in general, if we
 assume moreover that the dynamics is modified by the addition of a small
 random noise. But this is still an open problem to investigate.

 \bigskip

 \noindent {\bf Acknowledgements}\\
 This work was supported by the contract ANR MICA (2006-2009).

 \bigskip

 \noindent {\bf References}


\begin{thebibliography}{99}
 \bibitem{BC} 
 {Bulatov V V and Cai W},
 {\it Oxford University Press}, (2006).

 \bibitem{FTRC}
 {Fivel M, Tabourot L, Rauch E and Canova G R},
 {\it J. Phys.} IV, {\bf 8} (1998), 151-158.

\bibitem{HDK} 
{Hoc T, Devincre B and Kubin L P}, 
{\it In C. et al. Gundlach, editor, Riso National Laboratory, Denmark} (2004), 43-59. 

 \bibitem{K} 
 {Kr\"{o}ner E},
 {\it Erg. Angew. Math.} {\bf 5} (1958), 1-179, Berlin: Springer.

 \bibitem{K2} 
 {Kr\"{o}ner E},
 {\it Int. J. Solids and Structures} {\bf 38} (2001), 1115-1134.

 \bibitem{GB} 
 {Groma I and Balogh P},
 {\it Mat. Sci. Eng.} A {\bf 234-236} (1997), 249-252.

\bibitem{GCZ}
 {Groma I, Cikor F F and Zaiser M}, 
{\it Acta Mater.} {\bf 51} (2003), 1271-1281.


\bibitem{FE}
{El Hajj A and Forcadel N},
{\it Math. Comp.} {\bf 77} (2008), 789-812.

\bibitem{IJM}
{Ibrahim H, Jazar M and Monneau R},  
{\it C. R. Acad. Sci. Paris}, Ser I {\bf 346} (2008) 945-950.


 \bibitem{HZG}
 {Hochrainer T, Zaiser M and Gumbsch P},
 {\it Philosophical Magazine} {\bf 87} (8 \& 9) (2007), 1261-1282.

 \bibitem{M}
 {R. Monneau},
 {\it Interfaces Free Bound.}  {\bf 9} (2007), 383-409.

 \bibitem{FIM} 
 {Forcadel N, Imbert C and Monneau}, 
 {\it Discrete Contin. Dyn. Syst.} A {\bf 23} (3), to appear (March
 2009), and HAL: hal-00140545 (12-27-2007).

\bibitem{CIL}
 {Crandall M G, Ishii H and Lions P -L},
 {\it Bull. Amer. Math. Soc.} (N.S.) {\bf 27} (1992), 1-67.

\bibitem{AHLM}
{Alvarez O, Hoch P, Le Bouar Y and Monneau R},
{\it Arch. Ration. Mech. Anal.} {\bf 181} (3) (2006), 449-504.


 \end{thebibliography}
\end{document}